\documentclass[letterpaper, 11pt]{article}
\usepackage{amsmath, amsthm, amssymb}
\usepackage{graphicx}
\usepackage{xcolor}
\usepackage{algorithm}
\usepackage{algpseudocode}
\usepackage[colorlinks,linkcolor=blue,citecolor=blue,urlcolor=magenta,linktocpage,plainpages=false]{hyperref}
\usepackage{caption}
\usepackage{subcaption}
\usepackage{setspace}
\usepackage[left=0.9in, right=0.9in, top=0.9in, bottom=0.9in]{geometry}
\usepackage{authblk}
\usepackage{changepage}
\usepackage{enumitem}
\usepackage{setspace}

\usepackage{tikz}       % to use tikz drawing
\usetikzlibrary{calc}
\usetikzlibrary{plotmarks}
\usetikzlibrary{backgrounds}
\usetikzlibrary{arrows.meta,positioning}
\usetikzlibrary{decorations.pathreplacing}
\usetikzlibrary{patterns}
\usepackage{xcolor}
\usepackage{pgfplots}
\pgfplotsset{compat=newest}
\usetikzlibrary{positioning}

\usepackage{booktabs}
\usepackage[normalem]{ulem}

\newtheorem{remark}{Remark}

\definecolor{darkgreen}{RGB}{0, 128, 0}
\definecolor{orange}{RGB}{252, 138, 0}

\title{Extreme Strong Branching for QCQPs}
\date{}
\author[1]{Santanu S. Dey\thanks{\href{mailto:santanu.dey@isye.gatech.edu}{santanu.dey@isye.gatech.edu}}}
\author[2]{Dahye Han\thanks{\href{mailto:dahye.han@gatech.edu}{dahye.han@gatech.edu}}}
\author[3]{Yang Wang\thanks{\href{mailto:yang.wang@ce.gatech.edu}{yang.wang@ce.gatech.edu}}}

\affil[1,2]{School of Industrial and Systems Engineering, Georgia Institute of Technology}
\affil[3]{School of Civil and Environmental Engineering, Georgia Institute of Technology}

\begin{document}
% \onehalfspacing
\maketitle

\begin{abstract}
For mixed-integer programs (MIPs), strong branching is a highly effective variable selection method to reduce the number of nodes in the branch-and-bound algorithm. Extending it to nonlinear problems is conceptually simple but practically limited. Branching on a binary variable fixes the variable to 0 or 1, whereas branching on a continuous variable requires an additional decision to choose a branching point. Previous extensions of strong branching predefine this point and then solve $2n$ relaxations where $n$ is the number of candidate variables to branch. We propose \textit{extreme strong branching}, which evaluates multiple branching points per variable and jointly selects both the branching variable and point based on the objective value improvement. This approach resembles the success of strong branching for MIPs while additionally exploiting bound tightening as a byproduct. For certain types of quadratically constrained quadratic programs (QCQPs), computational experiments show that the extreme strong branching rule outperforms existing commercial solvers.
\end{abstract}

\section{Introduction}
The strong branching technique has been recognized as one of the most effective variable selection rules in the branch-and-bound algorithm for solving mixed-integer programs (MIPs) since its introduction \cite{applegate1995finding}; see also~\cite{achterberg2005branching,dey2024theoretical}. The method evaluates every fractional variable as a candidate variable to branch by fixing the candidate variable to either 0 or 1, solving relaxations at two child nodes, and evaluating the corresponding improvements in the objective function. The final branching variable is selected based on a composite score that combines the improvements from both branched child nodes, often measured by the product of the two gains in the child nodes.

Our work investigates whether the success of strong branching can be extended beyond binary variables. Although the extension to continuous variables is conceptually straightforward, its practical applicability remains limited and challenging. Unlike the case of branching on binary variables in MIP, branching on a continuous variable, known as \textit{spatial branching} \cite{smith1999symbolic}, requires an additional decision on the branching point. 

Recent advances in spatial branch-and-bound have introduced algorithmic strategies to address the fundamental challenges of variable and branching point selections. Violation transfer rule was proposed in \cite{tawarmalani2004global}, which systematically addresses constraint violations through strategic variable selection. The author of \cite{linderoth2005simplicial} developed a geometric approach that partitions the feasible region into triangles and rectangles, enabling finer spatial decomposition and convexification. Authors of~\cite{chen2017spatial} designed a specialized branching strategy for quadratically constrained quadratic programs (QCQPs) involving complex variables, addressing the unique difficulties in complex optimization domains. For bilinear problems, a branching rule that balances the violations across two child nodes was suggested in \cite{fischetti2020branch}, while specifically for bilinear bipartite programs (BBPs), authors of \cite{dey2019new} designed a branching rule that leverages second-order cone programming (SOCP) relaxations to guide variable selection decisions. More recently, a branching rule tailored for nonconvex separable piecewise linear functions was proposed in \cite{hubner2025spatial}. 

In addition, several works have adapted the strong branching technique from MIP to spatial branch-and-bound. One such variant is introduced in \cite{belotti2009branching}, in which a branching point is defined as a convex combination of the relaxed solution and the midpoint of a variable's bound, followed by solving $2n$ relaxations where $n$ is the number of candidate variables to evaluate the improvement associated with the preselected branching point for each variable. In \cite{kannan2022strong}, the authors proposed a learning-based approach to estimate both the branching point and branching variable expected to yield the largest improvement.

Our proposed approach performs an exhaustive search for a branching point for each variable using binary search while leveraging the bound-tightening technique during this process. This approach has two main benefits:
\begin{enumerate}%[leftmargin= 0.2in]
    \item Bound tightening is performed as part of the branching process.
    Although this resembles the spirit of the feasibility-based bound tightening or optimality-based bound tightening implemented in modern global solvers \cite{belotti2009branching, zhang2024solving}, one key distinction is that the bound reduction here arises directly from evaluating candidate branching points and thus is a byproduct of the branching process rather than a separate preprocessing or postprocessing step. 
    \item Since the extreme strong branching rule considers improvement in the objective function value, it resembles closely the objective-driven rationale that makes strong branching effective in MIP. The binary search restricts the number of evaluated points, providing an efficient trade-off between the computational search effort and improvement evaluation.
\end{enumerate}

\section{Extreme strong branching algorithm}
We consider a general nonlinear programming (NLP) of the form $\min \{ c^\top x \ | \ x \in P\}$. For simplicity of presentation, we have considered a linear objective here, but the method holds for nonconvex objective functions as well.

Within a branch-and-bound tree, a subproblem at node $k$ is given by:
\begin{equation}
\begin{aligned}
    \min \ & c^\top x \\
    s.t. \ & x \in P_k
    % & \overline{lb} \leq x \leq \overline{ub}
\end{aligned}
\end{equation}
where $P_k$ denotes the feasible region defined by the branching decisions taken along the path to node $k$. To obtain a dual bound of the problem, we consider $R_k$, a convex relaxation of $P_k$. For a branching decision on variable $x_i$ at threshold $\alpha$, we denote by:
$$R_k(x_i \leq \alpha) \ \ \text{and} \ \ R_k(x_i \geq \alpha),$$
the convex relaxation of $P_k \cap \{x_i \leq \alpha\}$ and $P_k \cap \{x_i \geq \alpha\}$, respectively.

Our proposed branching algorithm is provided in Algorithm~\ref{algo:extreme_branch} with two subroutines. The main algorithm is based on binary search. During the search, we may obtain a bound reduction. Finally, to choose the best branching variable and branching point combination, we compute branching scores. Each element is detailed in the following sections.

\subsection{Binary search}
\label{section:binary_search}
To identify the optimal branching threshold $\alpha$ for variable $x_i$, we employ a binary search over the interval defined by the current variable bounds $[\overline{lb}_i, \overline{ub}_i]$ for a fixed number of iterations. First, let us focus on the left child node problem corresponding to solving a relaxed subproblem of the form:
\begin{align*}
    obj_L(\alpha) = \min \{c^\top x \ | \ x \in R_k(x_i \leq \alpha) \}.
\end{align*}
Since the feasible region $R_k(x_i \leq \alpha)$ expands as $\alpha$ increases, $obj_L(\alpha)$ is a monotonically non-increasing function of $\alpha$. The first iteration begins with $\alpha = (\overline{lb}_i + \overline{ub}_i)/2$. If the subproblem at this $\alpha$ is feasible and the optimal objective function value $obj_L(\alpha)$ does not exceed the tree's best incumbent objective value (i.e., the tree's upper bound), then $obj_L(\alpha)$ is recorded and $\alpha$ is a candidate for branching. The search continues to the left for a smaller $\alpha$. On the other hand, if the subproblem is infeasible or $obj_L(\alpha)$ is greater than the tree's upper bound, the search direction is reversed to the right for a larger $\alpha$. We note that in this case, the lower bound of $x_i$ can be updated as further discussed in Section~\ref{section:bound_reduction}. This binary search is done until the binary search iteration limit.

The same binary search procedure is then applied symmetrically to the right child node corresponding to solving a relaxed subproblem of the form $obj_R(\alpha) = \min \{c^\top x \ | \ x \in R_k(x_i \leq \alpha) \}$. We note that $obj_R(\alpha)$ is now a monotonically non-decreasing function in $\alpha$.

After completing the binary search on both the left and right sides, we have a set of candidate branching points corresponding to different $\alpha$ values evaluated in both left and right directions. If any $\alpha$ has been explored on only one side, then additional subproblems are solved so that both $obj_L(\alpha)$ and $obj_R(\alpha)$ are available for every candidate $\alpha$. This ensures a consistent basis for computing the strong branching scores across the potential branching points, as further discussed in Section~\ref{section:branching_score}.

% Define Input and Output keywords
\algrenewcommand{\algorithmicrequire}{\textbf{Input:}}
\algrenewcommand{\algorithmicensure}{\textbf{Output:}}
\newcommand{\Input}{\Require}
\newcommand{\Output}{\Ensure}
\begin{algorithm}[H]
\footnotesize
\caption{Extreme Strong Branching Rule at Node $k$}
\label{algo:extreme_branch}
\begin{algorithmic}[1]
\Input The best upper bound of the tree ($obj_{ub}$) and optimal objective function value to the relaxation problem at node $k$,  $obj_{p} = \min \{c^T x \mid x \in R_k\}$
\Output The branching variable $x_{i^*}$ and the branching point $\alpha ^*$
\State $\alpha^* := 0$; $i^* := 0$; $score^* = -\infty$;
\For{each $i \in \{1, \ldots, n\}$}
    \State $lb_L = lb_R := \overline{lb} _i$; $ub_L = ub_R = \overline{ub}_i$; $B_L = B_R := \{\}$; $O_L = O_R = \{ \}$ \Comment{$B_{L}$ and $B_R$ represent sets of candidate branching points obtained, and $O_{L}$ and $O_R$ represents dictionaries of branching point and optimal objective function value pairs from solving the left and right child node problems respectively}
    \For{each $iter \in \{1, \ldots, iter_{\max}\}$}        
        % \State $(lb_L, ub_L, \overline{ub}_i, B_L, O_L)$ := 
        \State \texttt{binary\_search($i$, $\leq$, $lb_L$, $ub_L$, $\overline{ub}_i$, $B_L$, $O_L$)}
        % \State $(ub_R, lb_R, \overline{lb}_i, B_R, O_R)$ := 
        \State \texttt{binary\_search($i$, $\geq$, $ub_R$, $lb_R$, $\overline{lb}_i$, $B_R$, $O_R$)}
    \EndFor
    % \For{each $\alpha \in B_L$}
    %     \State $obj_L(\alpha) = O_L(\alpha)$ \Comment{Recall the optimal objective function value of the left side subproblem}
    %     \State \texttt{branch\_score($i$, $\alpha$, $\geq$, $obj_L(\alpha)$)} \Comment{Compute the branching score by solving the right side subproblem}
    % \EndFor    
    % \For{each $\alpha \in B_R$}
    %     \State $obj_R(\alpha) = O_R(\alpha)$ \Comment{Recall the optimal objective function value of the right side subproblem}
    %     \State \texttt{branch\_score($i$, $\alpha$, $\leq$, $obj_R(\alpha)$)} \Comment{Compute the branching score by solving the left side subproblem} 
    % \EndFor
    \For{each $\alpha \in B_L \cup B_R$}
        \State \texttt{branch\_score($i$, $\alpha$)} \Comment{Compute the branching score for each candidate branching points}
    \EndFor    
\EndFor
\State \Return $(i^*, \alpha^*)$
\end{algorithmic}
\end{algorithm}

\subsection{Bound reduction}
\label{section:bound_reduction}
First, we note the following conditions that we can tighten variable bounds.
\begin{remark}
\label{remark:bound_tighten}
If the problem $\min \{ c^T x \mid x \in R_k(x_i \leq \alpha) \}$ is infeasible or has an optimal objective value greater than the current upper bound of the branch-and-bound tree, then the lower bound of $x_i$ is at least $\alpha$. \end{remark}
% \begin{proof} We consider two cases: 
% \begin{enumerate}
%     \item Suppose $R_k(x_i \leq \alpha) = \emptyset$. Since $P_k \cap \{x \mid x_i \leq \alpha\} \subseteq R_k(x_i \leq \alpha)$,  this implies that $P_k \cap \{x \mid x_i \leq \alpha\} = \emptyset$. Hence, the optimal solution to the problem can only be in $P_k \cap \{x \mid  x_i \geq \alpha\}$.
%     \item Suppose $\min \{ c^T x \mid x \in R_k(x_i \leq \alpha) \} > obj_{ub}$. Since the true optimal objective value of the problem $obj^* \leq obj_{ub}$, it follows $obj^* < \min \{ c^T x \mid x \in R_k(x_i \leq \alpha) \}$ and thus $R_k(x_i\leq \alpha)$ does not contain the true optimal solution. Hence, $P_k \cap \{x \mid  x_i \leq \alpha\}$ also does not contain the true optimal solution.
% \end{enumerate}    
% \end{proof}

An analogous statement holds for the opposite branch: if $\min \{ c^T x \mid x \in R_k(x_i \geq \alpha) \}$ is infeasible or has an optimal objective value greater than the current upper bound of the branch-and-bound tree, then the upper bound of $x_i$ is at most $\alpha$. 

These observations provide a natural search direction for selecting the branching point $\alpha$ in Section~\ref{section:binary_search}. If the conditions of Remark~\ref{remark:bound_tighten} hold for the left branch, then values greater than the current $\alpha$ may be considered, knowing that $x_i \geq \alpha$ must hold. If the conditions of Remark~\ref{remark:bound_tighten} are not applicable, then values smaller than the current $\alpha$ may be considered for a potential bound tightening. 

% --- Subroutine: Binary Search ---
\begin{algorithm}[tbh!]
\footnotesize
\caption{\textbf{function} \texttt{binary\_search($i$, $\diamond$, $p_1$, $p_2$, $\overline{b}_i$, $B$, $O$)}}
\begin{algorithmic}[1]
\State $\alpha := (p_1 + p_2)/2$ \Comment{Compute the midpoint of the current interval as the next branching point}
\State $\text{obj} := \min \{ c^T x \mid x \in R_k(x_i \diamond \alpha) \}$ or $\infty$ if infeasible %\Comment{Solve the relaxation problem depending on the sign $\diamond$}
\If{$obj > obj_{ub}$} \Comment{If the subproblem's objective value exceeds the current best upper bound}
    \State $p_1 := \alpha$ \Comment{Shift the search direction towards the other side of the feasible region per Section~\ref{section:binary_search}}
    \State $ \overline{b}_i  := \alpha $ \Comment{Update the corresponding bound on variable $x_i$ per Section~\ref{section:bound_reduction}}
\Else
    \State $p_2 := \alpha$ \Comment{Otherwise, continue the search direction towards the same side of the feasible region}
    \State $O(\alpha) := obj$ \Comment{Record the optimal objective function value of this $\alpha$}
    \State $B := B \cup \{\alpha\}$ \Comment{Add $\alpha$ to potential branching point}
\EndIf
\end{algorithmic}
\end{algorithm}

\subsection{Branching score}
\label{section:branching_score} 
For each candidate threshold $\alpha$ generated during binary search, we compute a branching score to assess its effectiveness. Since we have explored different $\alpha$'s for the left and right child nodes, we might need to solve additional optimization problems to know the optimal objective function values for all $\alpha$'s that were considered. The score is defined as:
\[
score = \max \{ obj_L - obj_p, \epsilon \} \cdot \max \{ obj_R - obj_p, \epsilon \},
\]
where $obj_p$ is the optimal objective function value of the relaxation at node $k$ and $\epsilon > 0$ is a small stability constant. The branching variable and the branching point pair $(i^*, \alpha^*)$ with the highest branching score are selected as the branching decision at node $k$.

% % --- Subroutine: Branching Score ---
% \begin{algorithm}[tbh!]
% \footnotesize
% \caption{\textbf{function} \texttt{branch\_score($i$, $\alpha$, $\diamond$, $obj_{opp}$)}}
% \label{algo:branch_score}
% \begin{algorithmic}[1]
% \State $obj := \min \{ c^T x \mid x \in R_k(x_i \diamond \alpha) \}$ \Comment{Solve the relaxed subproblem for branching side $x_i \diamond \alpha$}
% \State $score_{new} := \max\{obj - obj_p, \epsilon\} \cdot \max \{obj_{opp} - obj_p, \epsilon\}$ \Comment{Compute the branching score as the product of bound improvements from both sides}
% \If{$score_{new} > score$} \Comment{If this branching variable and point combination gives better score}
% \State \textbf{then} $i^* = i$ and $\alpha^* := \alpha$ \Comment{Record the variable and branching point}
% \EndIf
% \end{algorithmic}
% \end{algorithm}

% --- Subroutine: Branching Score ---
\begin{algorithm}[tbh!]
\footnotesize
\caption{\textbf{function} \texttt{branch\_score($i$, $\alpha$)}}
\label{algo:branch_score}
\begin{algorithmic}[1]
\State \algorithmicif\ $\alpha$ is a key in $O_L$ \algorithmicthen\ $obj_L = O_L(\alpha)$ \algorithmicelse\ $obj_L = \min \{ c^T x \mid x \in R_k(x_i \leq \alpha) \}$
\State \algorithmicif\ $\alpha$ is a key in $O_R$ \algorithmicthen\ $obj_R = O_R(\alpha)$ \algorithmicelse\ $obj_R = \min \{ c^T x \mid x \in R_k(x_i \geq \alpha) \}$
\State $score_{new} := \max\{obj_L - obj_p, \epsilon\} \cdot \max \{obj_R - obj_p, \epsilon\}$ 
\State \algorithmicif\ $score_{new} > score^*$ \algorithmicthen\
$i^* = i$, $\alpha^* := \alpha$, and  $score^* = score_{new}$
\end{algorithmic}
\end{algorithm}

\begin{figure}[tbh!]
    \centering
    \begin{tikzpicture}
        \begin{axis}[
            width=13cm,
            height=6.3cm,
            xlabel={Branching Point ($\alpha$)},
            ylabel={Objective Function Value},
            xmin=0, xmax=1,
            ymin=0, ymax=9,
            axis line style={thick},
            label style={font=\bfseries\footnotesize},
            tick label style={font=\footnotesize},
            legend style={at={(1.15,0.3)}, anchor=north, legend columns=1, font=\footnotesize, draw=none},
            yticklabels={},
            xtick={0,0.25,0.31,0.375,0.5,0.56,0.625,0.75,1.0},
            grid=none,
            axis x line=bottom,
            xtick pos=bottom,
            ytick=\empty,
            y axis line style={draw=none},
            cycle list={},  % disable automatic color cycling
        ]

        % shaded gray areas (no legend)
        \addplot[draw=none, fill=gray!20, forget plot] coordinates {(0,0.01) (0.25,0.01) (0.25,9) (0,9)};
        \addplot[draw=none, fill=gray!20, forget plot] coordinates {(0.625,0.01) (1.05,0.01) (1.05,9) (0.625,9)};

        % obj_L line
        \addplot[red, thick, mark=*, mark options={fill=red,draw=red}, forget plot] 
            coordinates {(1.0,0) (0.5625,2) (0.5,2.5) (0.375,4) (0.3125,5.5) (0.25,7.8)}; % (0.625,1.5)

        % % obj_L line
        % \addplot[red, thick, mark=square, mark options={fill=red,draw=red}, forget plot] 
        %     coordinates {(0.5625,2) }; % (0.625,1.5)

        % obj_R line
        \addplot[blue, thick, mark=*, mark options={fill=blue,draw=blue}, forget plot] 
            coordinates {(0.0,0) (0.3125,1) (0.375,1.5) (0.5,3) (0.5625,4.5) (0.625,6.2) (0.75,7.5)}; %(0.25,0.3) 
    
        % horizontal dashed line
        \addplot[black, thick, dashed, forget plot] coordinates {(0,6) (1,6)};

        % ---- explicit legend images ----
        \addlegendimage{black, thick, dashed} % obj_UB
        \addlegendentry{$obj_{UB}$}

        \addlegendimage{red, thick, mark=*, mark options={fill=red,draw=red}, mark size=2pt} % obj_L
        \addlegendentry{$obj_L(\alpha)$}

        \addlegendimage{blue, thick, mark=*, mark options={fill=blue,draw=blue}, mark size=2pt} % obj_R
        \addlegendentry{$obj_R(\alpha)$}

        % iteration labels
        \node[blue] at (axis cs:0.48,3.65) {\footnotesize iter R1};
        \node[blue] at (axis cs:0.75,7.9) {\footnotesize iter R2};
        \node[blue] at (axis cs:0.6,6.7) {\footnotesize iter R3};
        \node[blue] at (axis cs:0.52,4.9) {\footnotesize iter R4};
        \node[red] at (axis cs:0.48,2.1) {\footnotesize iter L1};
        \node[red] at (axis cs:0.21,7.4) {\footnotesize iter L2};
        \node[red] at (axis cs:0.34,3.6) {\footnotesize iter L3};
        \node[red] at (axis cs:0.28,5.1) {\footnotesize iter L4};
        \node[draw=red, fill=none, circle, inner sep=4pt, thick, dashed]
          at (axis cs:0.5625,2) {};
        \node[draw=blue, fill=none, circle, inner sep=4pt, thick, dashed]
          at (axis cs:0.3125,1) {};
        \node[draw=blue, fill=none, circle, inner sep=4pt, thick, dashed]
          at (axis cs:0.375,1.5) {};

        \end{axis}
    \end{tikzpicture}
    \caption{\small Illustrative example of extreme strong branching for a single variable $x_i$. The x-axis represents candidate branching points $\alpha \in [0,1]$, and the y-axis shows the objective function value of the relaxed subproblem when branching at that point. The horizontal dashed line is the current best upper bound of the branch-and-bound tree, $obj_{UB}$. The procedure begins by exploring the left side $(x_i \leq \alpha)$. At iter L1, the left side problem is solved at the midpoint $\alpha = 0.5$. Because the problem is feasible and $obj_L(0.5) < obj_{UB}$, the algorithm continues to the left for a smaller $\alpha$ value. Hence, at iter L2, the left side is solved for the midpoint between 0 and the previous $\alpha$ value 0.5 is solved. We note that $obj_L(0.25) > obj_{UB}$, which implies $x_i \geq 0.25$. Hence, the left shaded area updates the lower bound of $x_i$. Per the algorithm, the direction is switched so we explore a larger value of $\alpha$. At iter L3 with $\alpha = (0.25 + 0.5)/2 = 0.375$, $obj_{L}(0.375) < obj_{UB}$ again, so at iter L4, we explore a smaller value of $\alpha = (0.25 + 0.375)/2 = 0.3125$. Once the left-side search terminates, the algorithm proceeds to the right side $(x_i \geq \alpha)$. The first right side iteration also starts at $\alpha = 0.5$ and subsequently moves towards larger $\alpha$ values whenever the right side problem remains feasible and below the upper bound. Otherwise, we evaluate smaller $\alpha$ values and tighten the upper bound of the variable. Shaded regions indicate these tightened bounds. After completing both left and right side binary searches, we have a set of candidate branching points. Finally, we solve the left or right child problems for these candidate points if they have not been solved before. These are marked by dashed lines around the dot. %We choose the branching variable and branching point combination that gives the highest improvement on both left and right child nodes.
}
\label{fig:branching}
\vspace{-0.3cm}    
\end{figure}
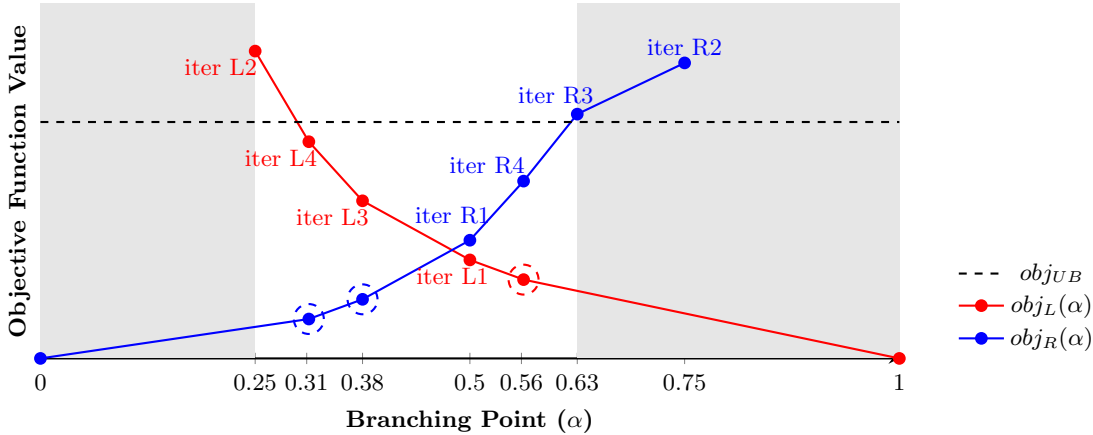

\paragraph{\textbf{Illustrative example}}
We illustrate the extreme strong branching rule applied for a single variable $x_i$ with bounds $[0,1]$. Figure~\ref{fig:branching} visualizes iteration steps taken in Algorithm~\ref{algo:extreme_branch}. The algorithm iteratively evaluates candidate branching points for the left child node problem $(x_i \leq \alpha)$ using binary search. Whenever points resulting in subproblems exceed the current best upper bound $obj_{UB}$, bounds of $x_i$ are updated, and the remaining feasible region is considered for deciding a new candidate branching point. The same technique is applied for the right child node problem $(x_i \geq \alpha)$. Finally, branching scores are calculated to select the branching variable and branching point.

\section{Branch-and-bound implementation}
We developed a custom branch-and-bound framework, building upon the Julia package \texttt{BranchAndBound.jl} \cite{bnbjulia}. 
Our implementation is designed to solve a general QCQPs of the form: 
\begin{equation}
\begin{aligned}
\label{eq:general_qcqp}
    \min_x \ & x^\top Q_{0} x + p_{0}^\top x \\
    s.t. \ & x^\top Q_k x + p_k^\top x \leq r_k, && \forall k \in [m] \\
    & lb \leq x \leq ub
\end{aligned}
\end{equation}
where $Q_0, Q_k \in \mathbb{R}^{n \times n}$, $p_0,p_k \in \mathbb{R}^n$, and $lb, ub$ are lower and upper bounds on the variable $x \in \mathbb{R}^n$.

% \paragraph{Upper bound initialization.} 
To obtain a high-quality feasible solution, we first run the global solver BARON~\cite{sahinidis1996baron} (version 25.3.19) for 60 seconds and pass its best incumbent solution as the initial primal heuristic to our branch-and-bound. This practice is similar to the standard practice of providing the best-known heuristic to isolate the evaluation of the dual bound improvement \cite{fischetti2020branch}. 

% \paragraph{Relaxation at each node.} 
The dual bound at each node is obtained by constructing a McCormick relaxation of the QCQP \cite{mccormick1976computability}. Specifically, each bilinear term $x_ix_j$ is replaced with an auxiliary variable $w_{ij}$ and the following McCormick inequalities are added:
\begin{equation}
\begin{aligned}
\label{eq:mccormick_relaxation}
%    \min \ &  <Q_{0}, w> + p_{0}^\top x \\
%    s.t. \ & <Q_k,  w> + p_k^\top x \leq r_k && \forall k \in [m] \\
    & w_{ij} - lb_j x_i - lb_ix_j + lb_i lb_j \geq 0 && \forall i, j \in [n]\\ 
    & w_{ij} - ub_j x_i - lb_ix_j + lb_i ub_j \leq 0 && \forall i, j \in [n]\\ 
    & w_{ij} - lb_j x_i - ub_ix_j + ub_i lb_j \leq 0 && \forall i, j \in [n]\\ 
    & w_{ij} - ub_j x_i - ub_ix_j + ub_i ub_j \geq 0 && \forall i, j \in [n].
%    & lb \leq x \leq ub
\end{aligned}
\end{equation}
Hence, each node in the branch-and-bound tree builds a McCormick relaxation of \eqref{eq:general_qcqp} with different local bounds of $lb$ and $ub$.

% \paragraph{Node selection and branching.} 
We use the best-bound rule for node selection, that is, we select the node with the smallest objective function value for a minimization problem.
%At each iteration, we select a node with the worst dual bound (that is, lowest objective function value for the minimization problem). The branching variable and branching point are then determined according to the branching rule that we consider: for our proposed extreme strong branching, Algorithm~\ref{algo:extreme_branch} is applied while alterative rules are also employed for benchmarking in Section~\ref{section:comp_results}. For experiments in model updating problem in Section~\ref{section:model_update}, feasibility-based bound tightening cuts are added where at each iteration whenever possible \cite{ryoo1995global, ryoo1996branch}. 

Once a node has been branched on, its two child nodes are added to the candidate pool while the original parent node is no longer considered. To reduce the computational overhead, the upper bound problem is solved every 10 iterations. Finally, Ipopt~\cite{wachter2006implementation} (version 3.14.17) was used as a local NLP solver to solve \eqref{eq:general_qcqp} and HiGHS~\cite{huangfu2018parallelizing} (version 1.11.0) was used as an LP solver to solve the McCormick relaxation with \eqref{eq:mccormick_relaxation}. At each node, a new LP model was built without warm starting the solver with any solutions obtained from solving LPs at different nodes. The maximum iteration of the branch-and-bound algorithm was limited to 100,000.

\section{Computational results}
\label{section:comp_results}
We evaluate the effectiveness of the proposed extreme strong branching rule in two settings: (i) benchmark QCQP instances from the literature \cite{bussieck2003minlplib}, and (ii) an application-driven problem from structural engineering, namely the finite element model (FEM) updating problem \cite{schreiber2023finite, dey2024aggregation}. All branch-and-bound implementations were coded in the Julia language version 1.11, and experiments were executed on a personal MacBook Air equipped with an Apple M3 chip (8-core CPU) with 8 gigabytes of RAM. We report the following metrics for performance:
\begin{itemize}
    \item Remaining optimality gap: using the best incumbent objective function value ($z^*$), we compute the remaining optimality gap as $\frac{|z^* - z_{lb}|}{|z^*|} \ (\%)$ where $z_{lb}$ denotes the lower bound on the objective function value obtained by different methods. %Since it measures the proximity to certifying optimality, lower values indicate better performance. 
    \item Solution time in seconds: The total elapsed time is measured in wall-clock seconds. For any instance unsolved within the time limit, the maximum time (3600 seconds) allowed is reported.
    \item Number of nodes explored: When using the custom branch-and-bound algorithm, we record the number of nodes processed in the tree. Commercial solvers that also rely on branch-and-bound tree search provide this metric in their logs.
\end{itemize}

\subsection{Benchmark QCQP instances}
% \textcolor{red}{say something about preliminary experiements showing the method works for bilinear problems...}
While the proposed rule was tested on a diverse set of QCQP instances, it proved especially effective for BBPs such as pooling problems and water management problems, which are the primary focus of our computational analysis. We considered problems with the number of variables at most 100. 
For these instances, to evaluate the effectiveness of the branching strategy, we compared three rules, where two rules were from previous literature on spatial branching. Both methods first consider the current infeasible relaxed solution $x^*$. We have not solved these with commercial solvers like Gurobi, as such solvers are highly optimized for standard instances.
\begin{itemize}
    \item \texttt{basic}: This is a basic version of the spatial strong branching. For each candidate branching variable $x_{i^*}$, a branching point is considered as a weighted combination of the midpoint of a variable's bound and the current relaxed solution. A branching score is computed following Algorithm~\ref{algo:branch_score} and a variable with the highest branching score is selected to branch \cite{belotti2009branching}.
    \item \texttt{balance}: The branching point and the branching variable pair are selected based on the rule specifically for bilinear problems \cite{fischetti2020branch}.
    \item \texttt{esb}: The branching variable and branching point are selected according to our proposed extreme strong branching rule in Algorithm~\ref{algo:extreme_branch}.
\end{itemize}

Table~\ref{table:minlp_instances} summarizes instances for which all three branching rules solved to $0.1\%$-optimality within the time limit. There are 35 instances of these and are grouped into four categories according to their solution time range. We report the number of instances solved (\# opt) as well as the arithmetic and geometric means of the solution times ($T_{ari}$ and $T_{geo}$, respectively). Overall, the \texttt{esb} rule is the only branching strategy that solves all of these instances under 100 seconds. While the \texttt{balance} rule performs competitively on some of the easier instances, we note that it requires solution time above 100 or even 1000 seconds for several of the more challenging instances. Both the arithmetic and geometric averages on the solution time of all instances also verify that \texttt{esb} rule consistently performs better. Table~\ref{table:unsolved_minlp_instances} provides the average remaining optimality gap for instances that were unsolved by at least one of the branching rules. For these more difficult instances, \texttt{esb} rule also provides the smallest remaining gap.

Detailed results are provided in Table~\ref{table_app:minlp_instances} in Appendix~\ref{section:appendix_computational_results}, where we further note that the \texttt{esb} rule explores drastically fewer nodes on average, which may not benefit much for problems requiring a relatively small-sized tree, but can outweigh the computational overhead for problems requiring a larger-sized tree by other rules. %Even when accounting for the fact that the \texttt{esb} rule solves roughly four times as many LP relaxations as the \texttt{basic} rule, it still explores fewer nodes overall (142 versus 825 by arithmetic average).

\begin{table}[tbh!]
\centering
\footnotesize
\caption{\small Average solve times for solved MINLP instances}
\label{table:minlp_instances}
\begin{tabular}{lrrrrrrrrrrr}
\toprule
 & \multicolumn{3}{c}{\texttt{basic}} &  & \multicolumn{3}{c}{\texttt{balance}} &  & \multicolumn{3}{c}{\texttt{esb}} \\
\cmidrule{2-4} \cmidrule{6-8} \cmidrule{10-12}
Time Range & \# opt & $T_{ari}$ & $T_{geo}$ &  & \# opt & $T_{ari}$ & $T_{geo}$ &  & \# opt & $T_{ari}$ & $T_{geo}$ \\
\midrule
(0,\ 10{]} & 17 & 5.09 & 5.07 &  & 20 & 4.08 & 3.97 &  & 19 & 4.09 & 3.88 \\
(10,\ 100{]} & 11 & 48.16 & 43.06 &  & 11 & 21.20 & 18.54 &  & 16 & 22.05 & 19.39 \\
(100,\ 1000{]} & 6 & 413.75 & 327.06 &  & 3 & 163.40 & 143.37 &  & 0 & n/a & n/a \\
(1000,\ 3600{]} & 1 & 2099.35 & 2099.35 &  & 1 & 1221.00 & 1221.00 &  & 0 & n/a & n/a \\
All & 35 & 148.52 & 24.11 &  & 35 & 57.89 & 10.33 &  & 35 & 12.30 & 8.09 \\
\bottomrule
\end{tabular}
\vspace{-0.2cm}
\end{table}

\begin{table}[tbh!]
\centering
\footnotesize
\caption{\small Average remaining optimality gap for unsolved MINLP instances}
\label{table:unsolved_minlp_instances}
\begin{tabular}{llll}
\toprule
 % & \multicolumn{5}{c}{remaining optimality gap} \\
 % \cmidrule{2-2} \cmidrule{4-4} \cmidrule{6-6}
\# inst & \texttt{basic} & \texttt{balance} & \texttt{esb} \\
\midrule
9 & 48.11\% & 51.21\% & 8.81\% \\
\bottomrule
\end{tabular}
\vspace{-0.2cm}
\end{table}

\subsection{Model updating problem}
\label{section:model_update}
% For FEM instances, we compared the entire branch-and-bound scheme with the extreme strong branching rule against general-purpose commercial solvers. Specifically, we compared against the following approaches:
% \begin{itemize}
%     \item \texttt{Gurobi}: a direct solve with a commercial solver Gurobi (version 12.0.1) \cite{gurobi}
%     \item \texttt{esb}: a branch-and-bound algorithm with our extreme strong branching rule
% \end{itemize}
For FEM instances, we compared the entire branch-and-bound scheme with the extreme strong branching rule against directly solving it with a general-purpose commercial solver Gurobi~\cite{gurobi} (version 12.0.1). The absolute constraint violation was set to $10^{-9}$ for Gurobi. Other parameters were kept as the default parameters. For both methods, we set a time limit of 3600 seconds.

Table~\ref{table:num_optimality_avg_time} summarizes the number of instances solved to 0.1\%-optimality within the time limit (\# opt), along with the arithmetic and geometric means of the solution times. The performance of the branch-and-bound scheme with our proposed extreme strong branching rule solves a larger number of instances (8 out of 20), exceeding Gurobi (5 solved). 
Table~\ref{table:remaining_gap_num_nodes} summarizes the three metrics across all instances, including instances that reached the time limit of 3600 seconds. Both the average remaining optimality gap and the number of nodes further show that the extreme strong branching rule can reach a smaller optimality gap with a smaller tree size. 
The entirety of the detailed computational results is provided in Table~\ref{table_app:remaining_gap_num_nodes} in Appendix~\ref{section:appendix_computational_results}.
We have also tested our instances against BARON and Counne~\cite{belotti2009couenne}, but they did not perform better than the extreme strong branching and achieved only marginal improvements in the lower bounds with their default parameters. Although the commercial solvers have highly engineered and sophisticated internal strategies, these results show that a simple yet effective branching rule can lead to substantial performance gains without the use of any advanced convexification techniques. 

\begin{table}[tbh!]
\centering
\footnotesize
\caption{\small Instances solved to optimality and average time limits.}
\label{table:num_optimality_avg_time}
\begin{tabular}{ccclccc}
\toprule
\multicolumn{3}{c}{\texttt{Gurobi}} &  & \multicolumn{3}{c}{\texttt{esb}} \\
\midrule
\# opt & $T_{avg}$ & $T_{geo}$ & \multicolumn{1}{c}{} & \# opt & $T_{avg}$ & $T_{geo}$ \\
\cmidrule{1-3} \cmidrule{5-7}
5 & 1418.81 & 325.99 & & 8 & 1345.98 & 970.51 \\
\bottomrule
\end{tabular}
\vspace{-0.3cm}
\end{table}

\begin{table}[tbh!]
\centering
\footnotesize
\caption{\small Summary across all instances on three evaluation metrics.}
\label{table:remaining_gap_num_nodes}
\begin{tabular}{lrrlrrlrr}
\toprule
 & \multicolumn{2}{c}{Remaining optimality gap} &  & \multicolumn{2}{c}{Solution time (sec)} & & \multicolumn{2}{c}{Number of nodes solved} \\
 \cmidrule{2-3} \cmidrule{5-6} \cmidrule{8-9}
 & \multicolumn{1}{c}{\texttt{Gurobi}} & \multicolumn{1}{c}{\texttt{esb}} &  & \multicolumn{1}{c}{\texttt{Gurobi}} & \multicolumn{1}{c}{\texttt{esb}} &  & \multicolumn{1}{c}{\texttt{Gurobi}} & \multicolumn{1}{c}{\texttt{esb}}\\
\midrule
Arithmetic mean & 11.26\% & 7.94\% &  & 3054.70 & 2698.39 &  & 8607384 & 58 \\
Geometric mean & 3.20\% & 1.65\% &  & 1974.82 & 2130.99 &  & 4533969 & 44 \\
\bottomrule
\end{tabular}
\vspace{-0.5cm}
\end{table}

\section{Conclusion}
In this paper, we introduced a new spatial branching rule, namely \textit{extreme strong branching}. The method combines binary search with bound tightening and extends the objective-driven efficiency of strong branching from MIP to MINLP. While applicable to general MINLPs, our preliminary computational experiments show that the rule is particularly effective for bilinear bipartite problems, which are special subcases of QCQPs. 

For problems with a large number of variables, our approach can have limitations, since all continuous variables are considered as candidate branching variables and the overhead computation to solve relaxation problems may outweigh the benefits. As a future direction, we would like to consider developing a reliability-based variant of extreme strong branching. Understanding why the method works particularly well for BBP is also another stream of future research.

\section*{Acknowledgment}
The authors would like to gratefully acknowledge the support of grant number 2211343 from the NSF CMMI.
\bibliographystyle{plain}
\bibliography{reference}

\newpage

\appendix
\section{Computational results}
\label{section:appendix_computational_results}

\begin{table}[tbh!]
\centering
\footnotesize
\caption{\small Remaining optimality gap, solution time, and number of nodes solved for MINLP instances: the lowest value either by the remaining optimality gap or the solution time among the three methods is highlighted in bold; if another method reached within 10\% of the lowest value, such a value is also highlighted.}
\label{table_app:minlp_instances}
\begin{tabular}{lrrrrrrrrrrr}
\toprule
 & \multicolumn{3}{c}{Remaining optimality gap} &  & \multicolumn{3}{c}{Solution time (sec)} &  & \multicolumn{3}{c}{Number of nodes} \\
instance & \texttt{basic} & \texttt{balance} & \texttt{esb} &  & \texttt{basic} & \texttt{balance} & \texttt{esb} &  & \texttt{basic} & \texttt{balance} & \texttt{esb} \\
\midrule
pooling\_adhya1pq & 0.10\% & 0.10\% & 0.10\% &  & 40.40 & 39.06 & \textbf{12.96} &  & 201 & 2899 & 15 \\
pooling\_adhya1stp & 0.10\% & 0.10\% & 0.09\% &  & 85.30 & 101.82 & \textbf{23.65} &  & 179 & 3551 & 11 \\
pooling\_adhya1tp & 0.10\% & 0.10\% & 0.11\% &  & 355.77 & \textbf{15.29} & \textbf{16.57} &  & 1073 & 979 & 21 \\
pooling\_adhya2pq & 0.09\% & 0.10\% & 0.04\% &  & 38.91 & 16.76 & \textbf{13.87} &  & 183 & 901 & 17 \\
pooling\_adhya2stp & 0.09\% & 0.10\% & 0.08\% &  & 106.83 & 100.57 & \textbf{22.40} &  & 217 & 3375 & 11 \\
pooling\_adhya2tp & 4.54\% & 0.10\% & 0.02\% &  & 55.63 & 15.09 & \textbf{14.13} &  & 215 & 885 & 17 \\
pooling\_adhya3pq & 0.10\% & 0.10\% & 0.00\% &  & 90.10 & \textbf{10.68} & 33.81 &  & 159 & 265 & 17 \\
pooling\_adhya3stp & 0.10\% & 0.10\% & 0.00\% &  & 222.75 & \textbf{52.88} & 64.78 &  & 153 & 945 & 11 \\
pooling\_adhya3tp & 2.45\% & 0.10\% & 0.01\% &  & 262.23 & \textbf{13.60} & 29.02 &  & 229 & 405 & 13 \\
pooling\_adhya4pq & 0.07\% & 0.10\% & 0.01\% &  & 14.67 & \textbf{8.33} & 15.26 &  & 15 & 151 & 5 \\
pooling\_adhya4stp & 0.07\% & 0.09\% & 0.03\% &  & 27.25 & \textbf{24.15} & 35.36 &  & 15 & 493 & 5 \\
pooling\_adhya4tp & 0.06\% & 0.10\% & 0.08\% &  & 573.45 & 19.89 & \textbf{15.15} &  & 421 & 531 & 5 \\
pooling\_bental4pq & 0.00\% & 0.00\% & 0.00\% &  & 5.03 & 4.00 & \textbf{3.57} &  & 7 & 9 & 5 \\
pooling\_bental4stp & 0.00\% & 0.00\% & 0.00\% &  & 5.19 & \textbf{3.97} & \textbf{3.71} &  & 7 & 13 & 5 \\
pooling\_bental4tp & 0.00\% & 0.00\% & 0.00\% &  & 5.03 & \textbf{3.87} & \textbf{3.57} &  & 7 & 11 & 5 \\
pooling\_bental5pq & 0.00\% & 0.00\% & 0.00\% &  & 4.39 & \textbf{2.97} & \textbf{3.00} &  & 1 & 1 & 1 \\
pooling\_bental5tp & 0.00\% & 0.00\% & 0.00\% &  & 4.71 & \textbf{3.04} & \textbf{3.00} &  & 1 & 1 & 1 \\
pooling\_foulds2pq & 0.00\% & 0.00\% & 0.00\% &  & 4.44 & \textbf{2.90} & \textbf{3.00} &  & 1 & 1 & 1 \\
pooling\_foulds2stp & 0.00\% & 0.00\% & 0.00\% &  & 6.40 & \textbf{4.03} & 7.15 &  & 1 & 1 & 1 \\
pooling\_foulds2tp & 0.00\% & 0.00\% & 0.00\% &  & 4.43 & \textbf{2.87} & \textbf{2.85} &  & 1 & 1 & 1 \\
pooling\_haverly1pq & 0.00\% & 0.00\% & 0.00\% &  & 5.18 & 3.90 & \textbf{3.45} &  & 5 & 9 & 5 \\
pooling\_haverly1stp & 0.00\% & 0.00\% & 0.00\% &  & 5.20 & \textbf{4.00} & \textbf{3.73} &  & 5 & 11 & 5 \\
pooling\_haverly1tp & 0.00\% & 0.00\% & 0.00\% &  & 5.17 & \textbf{3.90} & \textbf{3.79} &  & 9 & 9 & 5 \\
pooling\_haverly2pq & 0.00\% & 0.00\% & 0.00\% &  & 5.13 & 4.46 & \textbf{3.61} &  & 5 & 11 & 5 \\
pooling\_haverly2stp & 0.00\% & 0.00\% & 0.00\% &  & 5.20 & \textbf{3.97} & \textbf{3.66} &  & 5 & 13 & 5 \\
pooling\_haverly2tp & 0.00\% & 0.00\% & 0.00\% &  & 5.09 & 4.23 & \textbf{3.60} &  & 9 & 15 & 5 \\
pooling\_haverly3pq & 0.00\% & 0.08\% & 0.00\% &  & 5.27 & 4.23 & \textbf{3.54} &  & 7 & 33 & 5 \\
pooling\_haverly3stp & 0.00\% & 0.00\% & 0.00\% &  & 5.18 & \textbf{3.99} & \textbf{3.72} &  & 7 & 15 & 5 \\
pooling\_haverly3tp & 0.10\% & 0.00\% & 0.00\% &  & 5.56 & 3.98 & \textbf{3.52} &  & 25 & 9 & 5 \\
pooling\_rt2pq & 0.07\% & 0.10\% & 0.04\% &  & 37.62 & \textbf{10.30} & 12.40 &  & 173 & 627 & 13 \\
pooling\_rt2stp & 0.10\% & 0.09\% & 0.02\% &  & 57.19 & \textbf{15.51} & \textbf{14.97} &  & 117 & 563 & 7 \\
pooling\_rt2tp & 0.06\% & 0.05\% & 0.00\% &  & 50.13 & \textbf{4.32} & 10.01 &  & 225 & 33 & 11 \\
wastewater02m1 & 0.09\% & 0.10\% & 0.01\% &  & 32.57 & 287.80 & \textbf{5.73} &  & 547 & 67287 & 15 \\
wastewater02m2 & 0.10\% & 0.09\% & 0.00\% &  & 2099.35 & \textbf{4.68} & 18.49 &  & 703 & 217 & 19 \\
wastewater04m1 & 0.09\% & 17.10\% & 0.03\% &  & 961.47 & 1221.00 & \textbf{9.49} &  & 799 & 112681 & 23 \\
wastewater04m2 & 0.09\% & 0.10\% & 0.06\% &  & 3600.00 & 84.83 & \textbf{45.41} &  & 1087 & 4990 & 23 \\
wastewater05m1 & 47.53\% & 56.51\% & \textbf{0.08\%} &  & 3600.00 & 3600.00 & 3600.00 &  & 1521 & 65391 & 491 \\
wastewater14m1 & 58.43\% & 56.82\% & \textbf{37.43\%} &  & 3600.00 & 3600.00 & 3600.00 &  & 105 & 31921 & 609 \\
wastewater15m1 & 49.74\% & 53.11\% & \textbf{38.75\%} &  & 3600.00 & 3600.00 & 3600.00 &  & 8999 & 77303 & 1387 \\
waterund01 & 2.35\% & 37.79\% & \textbf{0.11\%} &  & 3600.00 & 2742.10 & 3600.00 &  & 11547 & 104743 & 2223 \\
waterund08 & 100.00\% & 79.10\% & \textbf{2.26\%} &  & 3600.00 & 3600.00 & 3600.00 &  & 285 & 37501 & 27 \\
waterund11 & 57.10\% & 57.10\% & \textbf{0.10\%} &  & 3600.00 & 3600.00 & \textbf{960.04} &  & 683 & 57693 & 189 \\
waterund17 & 55.15\% & 63.75\% & \textbf{0.39\%} &  & 3601.95 & 3600.00 & 3600.00 &  & 3095 & 49323 & 389 \\
waterund18 & 62.63\% & 56.63\% & \textbf{0.15\%} &  & 3600.09 & 3600.00 & 3600.00 &  & 3237 & 64743 & 599 \\
\midrule
Arithmetic mean & 10.03\% & 10.90\% & \textbf{1.82\%} &  & 854.51 & 683.02 & \textbf{605.36} &  & 825 & 15695 & \textbf{142} \\
Geometric mean & 0.45\% & 0.47\% & \textbf{0.18\%} &  & 67.13.35 & 31.21 & \textbf{24.75} &  & 67 & 302 & \textbf{14} \\
\bottomrule
\end{tabular}
\end{table}

\begin{table}[tbh!]
\centering
\footnotesize
\caption{\small Remaining optimality gap, solution time, and number of nodes solved by each method: the lower value either by the remaining optimality gap or the solution time achieved by two methods is highlighted in bold.}
\label{table_app:remaining_gap_num_nodes}
\begin{tabular}{lrrlrrlrr}
\toprule
 & \multicolumn{2}{c}{Remaining optimality gap} &  & \multicolumn{2}{c}{Solution time (sec)} & & \multicolumn{2}{c}{Number of nodes solved} \\
 \cmidrule{2-3} \cmidrule{5-6} \cmidrule{8-9}
Instance & \multicolumn{1}{c}{\texttt{Gurobi}} & \multicolumn{1}{c}{\texttt{esb}} &  & \multicolumn{1}{c}{\texttt{Gurobi}} & \multicolumn{1}{c}{\texttt{esb}} &  & \multicolumn{1}{c}{\texttt{Gurobi}} & \multicolumn{1}{c}{\texttt{esb}}\\
\midrule
in\_19\_48\_1 & 30.13\% & \textbf{27.51\%} &  & 3600.00 & 3600.00 &  & 9289972 & 105 \\
in\_19\_48\_2 & \textbf{0.10\%} & 28.29\% &  & 2436.81 & 3600.00 &  & 6882539 & 105 \\
in\_19\_48\_3 & 24.15\% & \textbf{14.74\%} &  & 3600.00 & 3600.00 &  & 14495254 & 93 \\
in\_19\_48\_4 & 0.05\% & 0.06\% &  & \textbf{1.69} & 424.67 &  & 3147 & 13 \\
in\_19\_48\_5 & \textbf{0.98\%} & 2.10\% &  & 3600.00 & 3600.00 &  & 16777932 & 107 \\
in\_19\_48\_6 & 10.64\% & \textbf{0.10\%} &  & 3600.00 & 3009.85 &  & 12769840 & 119 \\
in\_19\_48\_7 & \textbf{5.59\%} & 16.30\% &  & 3600.00 & 3600.00 &  & 15071785 & 103 \\
in\_19\_48\_8 & 0.99\% & \textbf{0.88\%} &  & 3600.00 & 3600.00 &  & 15890286 & 71 \\
in\_19\_48\_9 & 9.15\% & \textbf{3.93\%} &  & 3600.00 & 3600.00 &  & 10464904 & 89 \\
in\_19\_48\_10 & 2.30\% & \textbf{0.00\%} &  & 3600.00 & 1002.23 &  & 13998031 & 17 \\
in\_21\_54\_1 & 35.28\% & \multicolumn{1}{r}{\textbf{16.33\%}} &  & 3600.00 & 3600.00 &  & 8509977 & 51 \\
in\_21\_54\_2 & 0.09\% & 0.00\% &  & \textbf{277.49} & 387.46 &  & 208340 & 13 \\
in\_21\_54\_3 & 38.00\% & \textbf{12.39\%} &  & 3600.00 & 3600.00 &  & 6459876 & 21 \\
in\_21\_54\_4 & 6.53\% & \textbf{2.20\%} &  & 3600.00 & 3600.00 &  & 5519562 & 63 \\
in\_21\_54\_5 & 38.24\% & \textbf{24.77\%} &  & 3600.00 & 3600.00 &  & 6942537 & 63 \\
in\_21\_54\_6 & 0.09\% & 0.00\% &  & 935.95 & \textbf{435.51} &  & 633461 & 15 \\
in\_21\_54\_7 & 15.88\% & \textbf{9.06\%} &  & 3600.00 & 3600.00 &  & 9509896 & 33 \\
in\_21\_54\_8 & 0.10\% & 0.01\% &  & 3442.11 & \textbf{675.41} &  & 4893953 & 11 \\
in\_21\_54\_9 & 5.38\% & \textbf{0.10\%} &  & 3600.00 & 1747.07 &  & 7447186 & 27 \\
in\_21\_54\_10 & 1.63\% & \textbf{0.07\%} &  & 3600.00 & 3085.60 &  & 6379203 & 39 \\
\midrule
Arithmetic mean & 11.26\% & \textbf{7.94\%} &  & 3054.70 & \textbf{2698.39} &  & 8607384 & \textbf{58} \\
Geometric mean & 3.20\% & \textbf{1.65\%} &  & 1974.82 & \textbf{2130.99} &  & 4533969 & \textbf{44} \\
\bottomrule
\end{tabular}
\end{table}

\end{document}